%% file: manuscript.tex
\documentclass[runningheads]{llncs}
\usepackage[T1]{fontenc}
\usepackage{graphicx}
\usepackage{hyperref}
\usepackage{color}

\usepackage[utf8]{inputenc}
\usepackage{xspace}
\usepackage{environ}
\usepackage{array}

\newif\ifshowdisc
\showdisctrue

\usepackage{xcolor}

\NewEnviron{discussion}{
\ifshowdisc
\hrule
START DISCUSSION
\BODY
END DISCUSSION
\hrule
\else
\fi
}

\input{notations}

\title{Knowledge management in House of Graphs}
\begin{document}
%
%
%

\author{Gauvain Devillez\inst{1}\orcidID{0000-0002-3931-1610} \and
Sven D'hondt\inst{2} \and
Jan Goedgebeur\inst{2,3}\orcidID{0000-0001-8984-2463}}
\authorrunning{G.\ Devillez et al.}
%
\institute{Computer Science Department, University of Mons, 7000 Mons, Belgium\\ 
\email{gauvain.devillez@umons.ac.be}\\ 
\and
Department of Computer Science, KU Leuven Kulak, 8500 Kortrijk, Belgium\\
\email{dhondtsven36@gmail.com,jan.goedgebeur@kuleuven.be}
\and
Department of Mathematics, Computer Science and Statistics, Ghent
University, 9000 Ghent, Belgium}

\maketitle              
\begin{abstract}
\Thog is an online database of graphs which can be accessed at \url{https://houseofgraphs.org/}. It serves as a central repository for complete lists of graphs for various graph classes. However, its main feature is a searchable database of so-called ``interesting'' graphs.
The development of the original House of Graphs started in 2010 and it was completely rebuilt in 2021-2022. 
Each graph in the database is accompanied by a significant amount of meta-data
such as a name, drawings, precomputed graph invariants, and comments. Given this
amount of information and the importance of reliability in the scientific world, robust data management is essential to ensure accuracy and
consistency across the database.
In this article, we therefore focus on knowledge management in \thog and describe the inner workings of \thog and how we ensure that its data is coherent, qualitative and stable.

\end{abstract}

\section{Introduction}

\Thog is, at its core, a community-driven database of undirected simple graphs (i.e., without semi-edges, parallel edges or loops). It is searchable
through its own website (\url{https://houseofgraphs.org/}) and can be augmented by registered users. As the number of graphs of a given order grows very fast, even if we limited ourselves to more restricted graph classes such as regular graphs\footnote{See \url{https://oeis.org/A001349} for the counts of (connected) graphs and \url{https://oeis.org/A005177} for the counts of (connected) regular graphs.}, \thog focuses on  ``interesting'' graphs, trading exhaustivity for curation. This raises
difficulties to ensure good knowledge management.

Most researchers will agree that not all graphs are equally interesting. Moreover, if a graph is ``interesting'' for one problem, it is often interesting for several -- sometimes seemingly unrelated -- problems. The Petersen graph is a prime example of this.

\Thog does not give a precise definition of which graphs are considered ``interesting''. Most graphs in \thog are extremal graphs for a given property or counterexamples to a conjecture. An important aspect of \thog is that users can upload their own interesting graphs. Basically, if a user finds a certain graph interesting enough to go through the effort of uploading the graph, we consider it as ``interesting'' and suitable for \thog.

For every new graph uploaded to \thog, several graph invariants are computed and stored in the database. Users can then perform queries to search for graphs with certain properties.

Next to the dynamic searchable database of ``interesting'' graphs, \thog also offers static exhaustive lists of all graphs of a given graph class (e.g.\ cographs, cubic graphs, quartic graphs, snarks, trees, etc.) up to the orders for which it was still feasible to store them. This is to accommodate users who are looking for exhaustive lists of graphs. They can download these lists and then process and filter them via their own tools. Previously, quite some static exhaustive lists of graphs already existed, but they were scattered over the internet and were typically hosted on the homepages of individual researchers. So another goal of \thog is to act as a central repository for static exhaustive lists of graphs for various graph classes.

The development on \thog started in 2010 and a first article describing \thog was published in~\cite{hog1}. In 2021-2022 \thog was completely rebuilt from scratch using more modern and future-proof frameworks as many of the underlying libraries of the original 2010 House of Graphs were no longer supported. Moreover, the old House of Graphs website did not render well on mobile devices. Next to that, the new House of Graphs was extended with a lot of new functionalities (e.g.\ searching for subgraphs, searching by invariant formulas, an improved graph drawing tool, several new graph invariants, etc.). All user-uploaded (meta-)data from the old House of Graphs was migrated to the new House of Graphs and the new website was described in the paper~\cite{hog2}.

The two aforementioned papers describe the functionalities of \thog from a non-technical (mathematical) user point of view. In this paper, however, we now describe the inner workings of \thog and how we ensure that its data is
coherent, qualitative and stable. Through this we intend to lift the veil and show users how \thog works below the surface and also share our experiences and best practices with other developers who are working with mathematical datasets. Indeed, as explained by Ber{\v{c}}i{\v{c}} in~\cite{katjaData}, the main challenges for the data quality of mathematical datasets concern correctness, completeness, consistency, and accessibility of the data. Examples of other large mathematical databases
include \textit{The On-Line Encyclopedia of Integer Sequences
(OEIS)}~\cite{oeis} and \textit{The {L}-functions and modular forms database
(LMFDB)}~\cite{lmfdb}. We refer the interested reader to
\textit{MathBases}~\cite{mathbases} for an index of mathematical databases.

The rest of this paper is organised as follows.
In Section~\ref{sec:hog} we first describe \thog
and its functionalities. Then, in Section~\ref{sec:knowledge} we describe how data is
entered in \thog. In Section~\ref{sec:quality} we go into detail about how the quality
and correctness of this data is ensured. Next, in Section~\ref{sec:exploration} we describe how the data can be searched and some
special measures we took to balance the users' ability to explore the database and system
robustness. In Section~\ref{sec:stability} we describe the design decisions we took to
preserve knowledge stability. Finally, we end in Section~\ref{sec:conclusion} with some concluding remarks.

\section{\Thog}\label{sec:hog}

\begin{figure}[h!]
  \centering
  \includegraphics[width=0.9\textwidth]{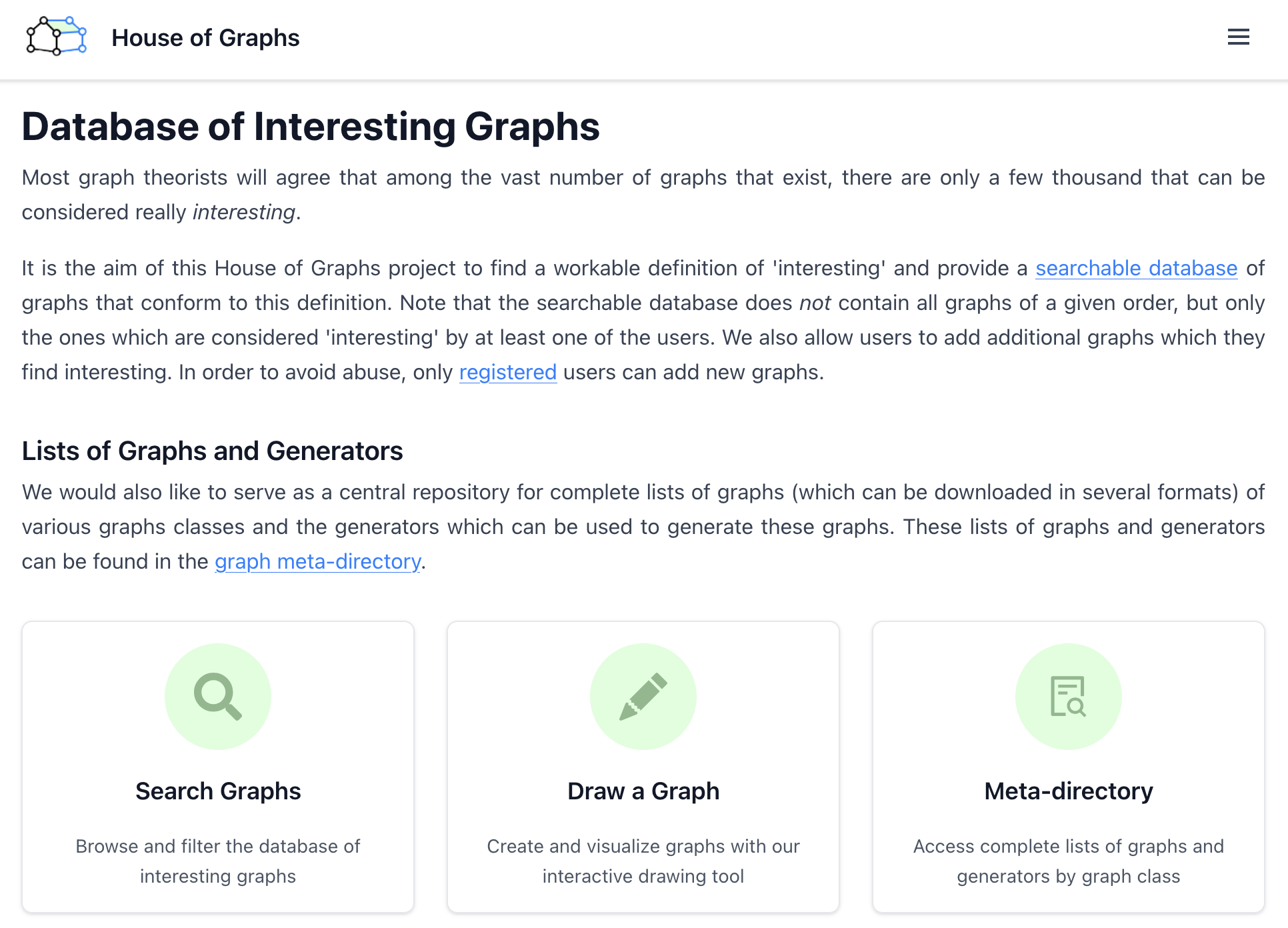}
  \caption{\Thog landing page.}
  \label{fig:landing}
\end{figure}

\Thog is a popular database of interesting graphs. At the time of writing, it has more than 600 registered users and many more unregistered visitors.
It provides multiple functionalities including
searching for graphs, submitting new graphs in multiple formats (including
drawing them), and complete lists of graphs via the \emph{meta-directory}. 
These options are visible on the landing page shown in Figure~\ref{fig:landing}.

\subsection{Interesting graphs}

By their combinatorial nature, the number of graphs of a given order grows exponentially fast. 
This is a nice illustration of the ``combinatorial explosion''. For example, there are already more than 165 billion pairwise non-isomorphic undirected graphs with 12 vertices\footnote{See \url{https://oeis.org/A000088}.}. Therefore, storing \textit{all} graphs is a herculean task and not feasible in practice.

\Thog's main idea to tackle this issue is to only store what we refer to as ``interesting'' graphs. 
As stated in the introduction, this concept takes many definitions depending on the person or the
context, so we intentionally do not give a formal definition. Informally, this means graphs that appear in the literature as
representative of a special family of graphs, counterexamples to 
conjectures, optimal structures with regard to some graph invariant, etc.

Each ``interesting'' graph is equipped with meta-data such as precomputed values for a wide selection of graph
invariants (at the time of writing, \thog offers 52 graph invariants), one or more drawings, and a unique numerical identifier. Next to that, users can
provide an optional name and additional drawings. The reason why the person who uploaded the graph (or other users) find the graph particularly interesting, is described through comments, which are displayed under the graph
information. On top of that, the user who submitted the graph can mark a selection of the available invariants as ``interesting'' for this graph (e.g.\ because the graph is extremal with regard to those invariants).

\begin{figure}[h!]
  \centering
  \includegraphics[width=0.9\textwidth]{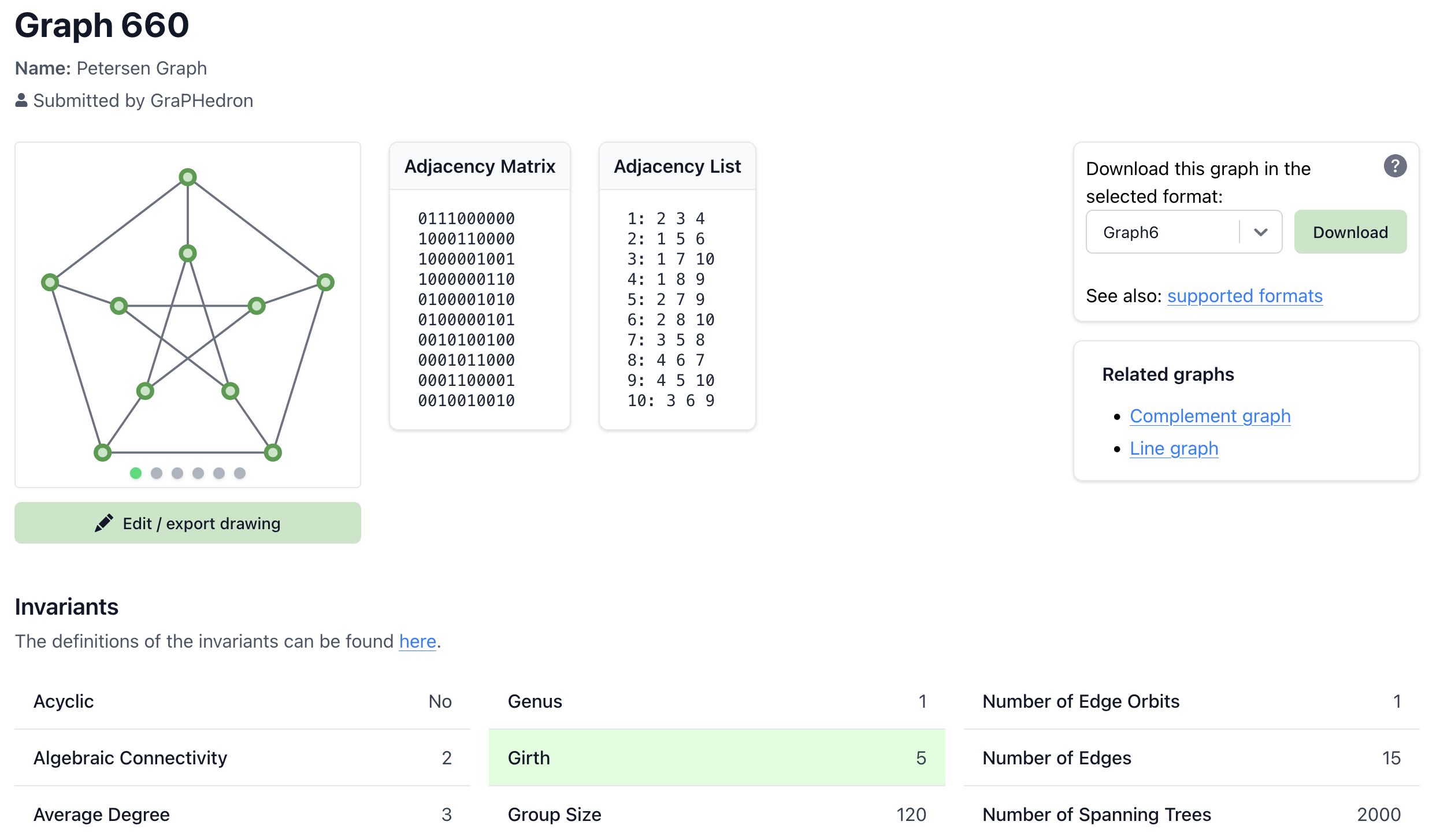}
  \caption{Excerpt of the detailed overview of a given graph. Invariants with a green background have been marked as particularly ``interesting'' for this graph by the person who uploaded the graph.}
  \label{fig:petersen}
\end{figure}

Figure~\ref{fig:petersen} shows the information provided for a given graph.
As can be seen, the graph can be
downloaded in multiple export formats. To preserve space, comments are not shown here (and only a selection of the graph invariants are shown).
This particular graph is named the \emph{Petersen Graph} and has the numerical House of Graphs identifier
660. If present in the database, \thog also provides links to related graphs
such as the complement or the line graph.

Note that \thog is currently restricted to simple \textit{undirected} graphs. There are several reasons for this. First of all, in graph theory most researchers focus on undirected graphs (although admittedly in computer science and other applications directed graphs are often used).
Furthermore, directed graphs suffer much more severely from a combinatorial explosion than undirected graphs (e.g.\ the number of tournaments alone -- constrained orientations of a complete graph -- already grows exponentially\footnote{See \url{https://oeis.org/A000568}.})
and technological challenges would arise should we try to store them. Moreover, the concept of ``interesting'' directed graphs is even less clear at the moment. In order for a directed graph to be interesting, the underlying undirected graph is just as important as the orientation of its edges. For example, there are many ``named'' undirected graphs in the literature but we are not aware of named directed graphs of graph-theoretical interest. Moreover, many invariants or invariant computers for undirected graphs do not work (or are not defined) for directed graphs.  However, this is certainly an interesting direction that we could keep in mind as possible future work. Though, if we do this, this will most likely result in a separate ``House of Directed Graphs''.

\subsection{Meta-directory}

\begin{figure}[h!]
  \centering
  \includegraphics[width=0.9\textwidth]{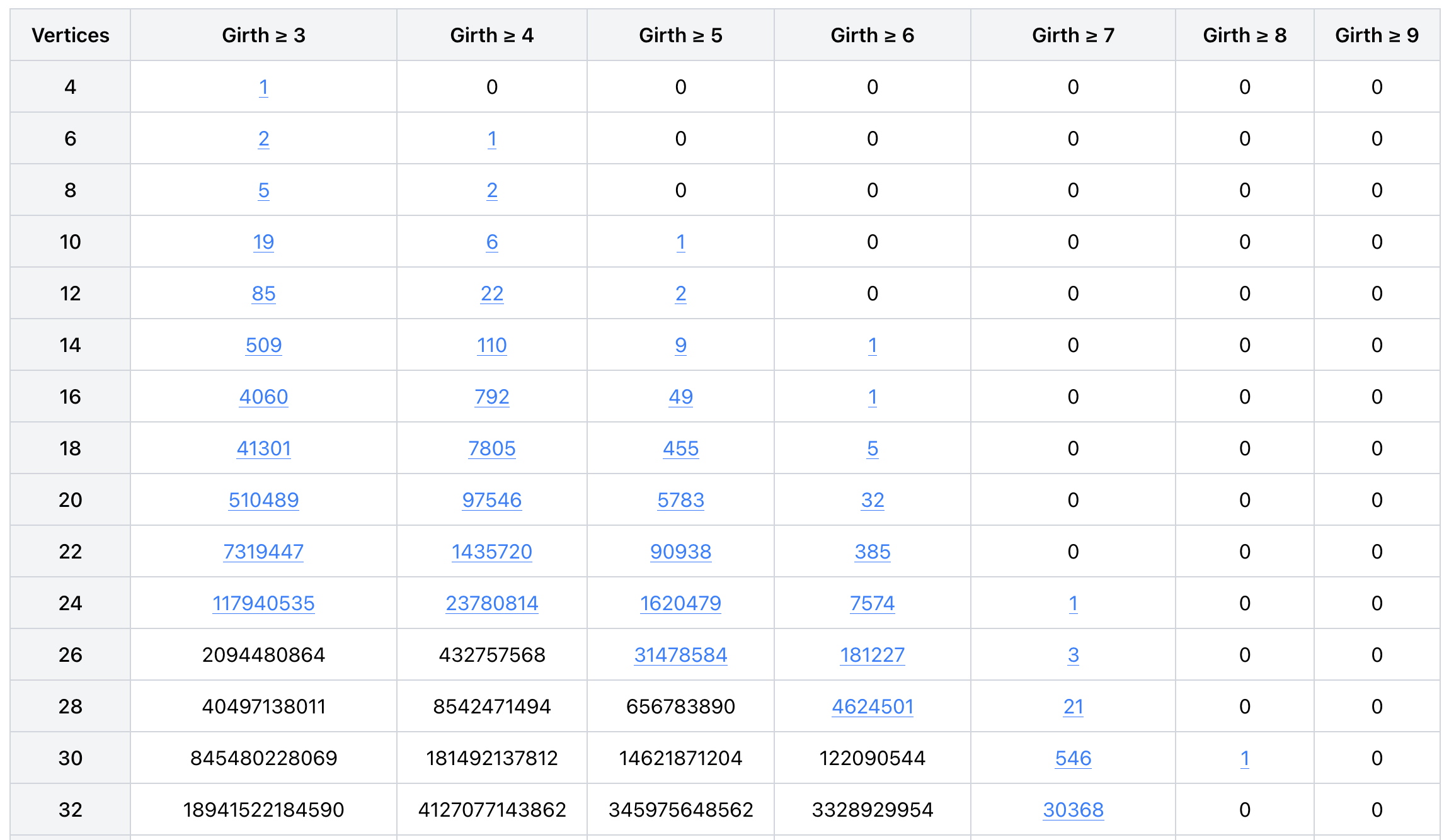}
  \caption{Sample of the cubic graphs page in the \textit{meta-directory}.}
  \label{fig:meta-cubic}
\end{figure}

Because the database is non-exhaustive by design, \thog also offers a
\emph{meta-directory} containing static exhaustive lists of graphs for several graph families as
well as links to the source code of the generators to generate them. Figure~\ref{fig:meta-cubic}
shows the meta-directory page dedicated to cubic graphs. For storage reasons, not \textit{all} lists are offered as a downloadable list. When the number of graphs in the list for a given order is small enough, a link
to this list is available (in blue in the figure). Otherwise, just the total number
is provided (in black). In some cases, however, the number of graphs of a given order is unknown (because it is computationally infeasible) and a question mark is shown. When available, \thog also
provides links to generators so the users can generate these graphs themselves (up to higher orders than what our storage allows), if they wish.

\subsection{Exploring the database}

When looking into ``interesting'' graphs, the focus of the user shifts from the static exhaustive collection of graphs of a given family from the \textit{meta-directory} to exploring individual graphs. 
\Thog therefore provides a search functionality for ``interesting'' graphs that allows
combining multiple types of constraints. We thereby aim to answer two needs of the users: 1) finding graphs that have specific properties, and 2) finding information about
a specific graph. Compared to an exhaustive dataset, the first goal might seem
less useful in a database that might not contain some graphs. However, an exhaustive dataset is likely to produce a very large number of graphs, often too many to handle. Having only ``interesting'' graphs will yield a smaller sample of promising graphs and allows to go up to higher orders while making it more likely to contain graphs of interest for the user's specific research problem.

The first set of available search constraints
focuses on invariant values. One can search for graphs by giving bounds on their
invariant values, exact values or parity. For boolean invariants, forming a class of
graphs (e.g.\ biparte graphs), one can ask for either inclusion or non-inclusion. A more advanced type
of constraint takes the form of a boolean formula concerning the values of numerical
invariants.

Further constraints relate to the information which was provided by users. One can search for graphs where a certain invariant was marked as ``interesting''. Comments and
names provided by users can also be searched. One can ask for some text to be
present either in one of the comments associated with the graphs or in their name.
A last constraint that is computed online is subgraph inclusion or exclusion.
That is, one can search for all graphs that contain or do not contain a given subgraph.

\begin{figure}[h!]
  \centering
  \includegraphics[width=0.9\textwidth]{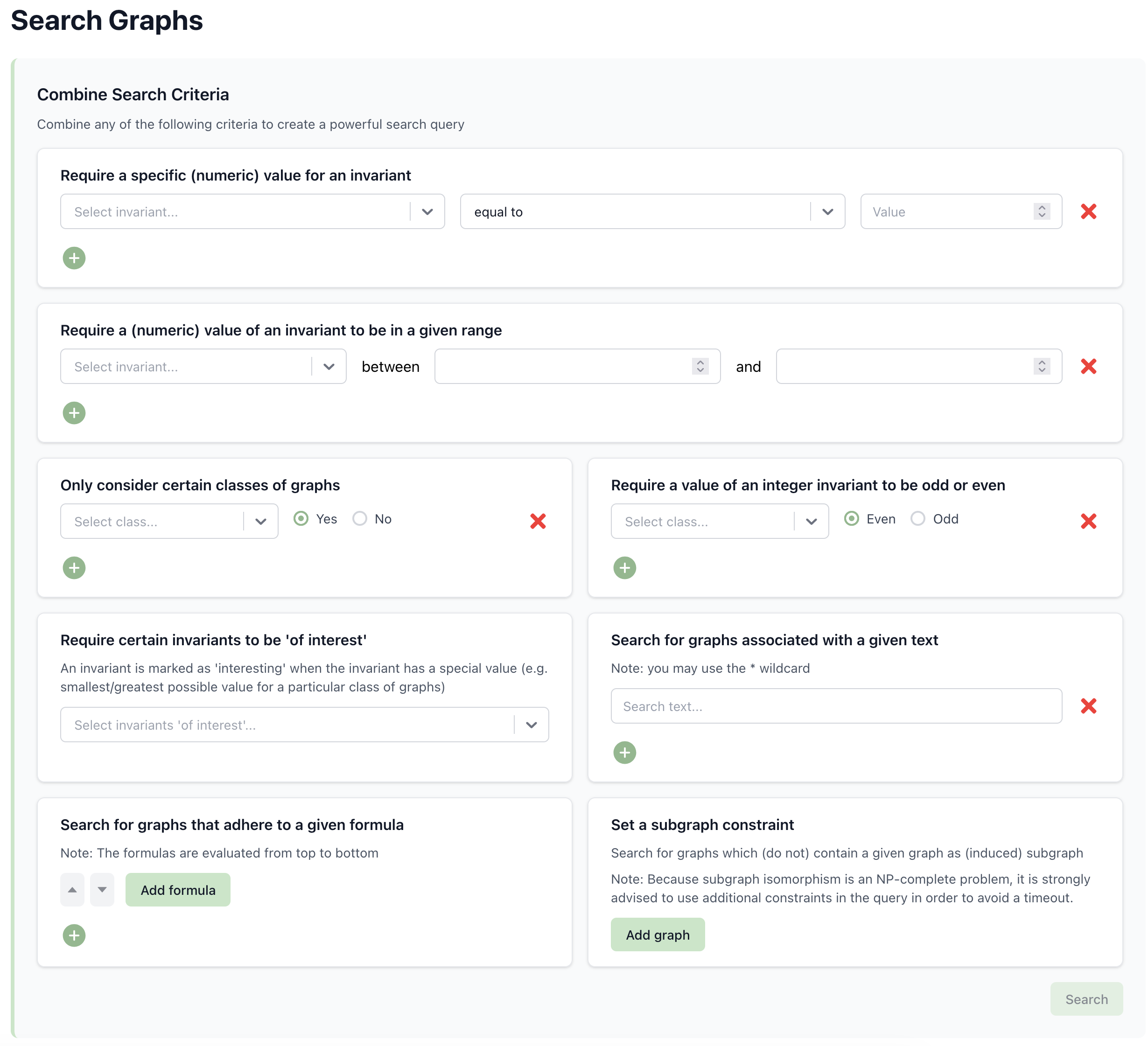}
  \caption{Overview of the search page.}
  \label{fig:search}
\end{figure}

All of these constraints can be combined, including also with multiple constraints of
each type. Figure~\ref{fig:search} displays the search page.

It is also possible to draw a graph and search if a graph isomorphic to this graph was already present in the database (e.g.\ to find out in what context other users came across this graph).

\section{Knowledge management}\label{sec:knowledge}

In this section we describe how knowledge management is handled in \thog. In Section~\ref{sec:initial-data}, we first describe how the database of ``interesting'' graphs was originally initialised and in Section~\ref{sec:comm-contr} how it is being augmented by the community.

\subsection{Initial data}\label{sec:initial-data}

Although the goal is to have a community-driven database that is continuously extended by users, it was important to initialise the database with a significant number of ``interesting'' graphs to have a viable starting point.

Initially, several well-known graphs from the literature were uploaded such as the Coxeter graph, Heawood graph, Petersen graph, etc. We amongst others uploaded all ``named'' graphs from \textit{Wolfram MathWorld}~\cite{mathworld} (which were kindly provided to us by Eric Weisstein) to \thog.
However, fully reviewing the extensive and ever-growing literature in graph theory seemed infeasible.

We therefore also used the system \textit{\gph} of M\'elot~\cite{melot2008facet}, now succeeded by \textit{PHOEG}~\cite{phoeg}, to obtain lists of graphs which are interesting for their
extremal properties. In particular, \textit{\gph} starts from the complete enumeration of graphs with
a fixed number of vertices. Each graph is then projected in the plane using the values of
two numerical invariants as coordinates. From this cloud of points, the convex hull is computed, thereby providing the extremal points and the graphs
associated with them (see~\cite{hog1} for more details). Note that this requires the values of the selected invariants for the complete enumeration. \textit{\gph} used graphs up to 10 vertices, that is, more than 12 million graphs.

Finally, we also uploaded several ``interesting'' graphs which we obtained in our own research, such as snarks which occurred as counterexamples~\cite{brinkmann2013generation} or extremal graphs for Ramsey-type problems~\cite{goedgebeur2013new}.

\subsection{Community contribution}\label{sec:comm-contr}

In order to extend (and keep extending) the list of interesting graphs, \thog was then opened for public submission in 2012, turning it into a community-driven effort. Thanks to the community, from an initial list of about 1\,500 graphs, \thog grew significantly and currently contains more than 29\,000 graphs at the time of writing. Despite this increase, we are far from even the number of graphs on say 10 vertices, which illustrates how rare ``interesting'' graphs are.

However, user-submitted content requires a careful design as we want to limit the restrictions imposed on users while
also preventing abuse or (unintentional) user errors. First of all, users need to register and create an account before they can add new graphs to the database (or add comments or new graph drawings of existing graphs).

It is important that users do not insert
thousands of randomly generated huge graphs, thereby defeating the purpose of
storing only ``interesting'' graphs. But the submission of individual graphs should
be straightforward.

To this end, users are only allowed to submit one graph at a time and their number of vertices must be between 1 and 250. (The \emph{null
graph}, with zero vertices, can also be seen as ``interesting'' on its own, but it also brings
 issues and it is therefore ignored in most of the literature~\cite{harary2006null}).

The upper limit has two reasons. Computationally speaking, as many graph invariants are known to be NP-hard to compute, it could take an extremely long time to compute these invariants for large graphs and this could also risk overloading the server. Second, very large graphs are generally very difficult to visualise without dedicated drawing algorithms.

Because we want that for every graph a reason is provided why it is interesting, the system ensures that, when a user uploads a new graph, they must also provide a comment with such a reason. If the graph would already be present in
\thog, the user is invited to add their own comment(s) about why they find this graph
interesting (if that reason was not mentioned yet in the previous comments).

The resulting dataset is relatively small (around 380\,MB in total), so its size is manageable in comparison with an exhaustive dataset of all graphs of a given order (e.g.\ storing all graphs of order 10 can already amount to 100\,MB). The median number of vertices is 29 with 75\% being under 45 vertices and the median number of edges is 68. In~\cite{hogISA} an \textit{instance space analysis} can be found with more insights about the properties of the graphs in \thog.
This shows that, at least in our context of ``interesting'' graphs, the main challenge is curation rather than storage.

\section{Quality assurance}\label{sec:quality}

Besides storing the knowledge, \thog must ensure that it is correct, coherent
and relevant to the graph theory community. Section~\ref{sec:corr-coher} describes how we ensure that the
data is correct and coherent with itself. In Section~\ref{sec:pert-inform}, we
describe the thought process to decide which information is included.

\subsection{Correctness and coherency}\label{sec:corr-coher}

A critical requirement is to ensure that the knowledge offered is correct and
coherent. That is, it does not contain errors and does not contradict itself.

To ensure coherency, it is important to avoid duplicates. For graphs, this
implies dealing with the graph isomorphism problem. \Thog relies on the
\emph{nauty} library~\cite{nauty} to compute a \emph{canonical labelling} of the
vertices of the graph.
That is, an ordering that is unique up to isomorphism (i.e.\ two graphs have the same canonical labelling if and only if they are isomorphic). Note that there are many options to define a canonical labelling. One commonly used option is to define this as the labelling which yields the smallest (or largest) binary number when you concatenate the rows of the adjacency matrix of the given labelled graph. However, we opted to use the canonical labelling of \emph{nauty} (which is rather technical to define), as it uses several clever heuristics so it can be computed very quickly in practice for most graphs (definitely for the orders we are considering in \thog).

For every graph in the database, a canonical labelling has been computed and the
canonically labelled graph is encoded and stored in the database in
\textit{graph6 format}~\cite{graph6}, i.e., a graph format which is often used
to encode graphs and is an encoding of its adjacency matrix compressed into a
sequence of printable characters.
When a new graph is submitted, its
canonical labelling is first computed using \emph{nauty}~\cite{nauty} and the
adjacency matrix of the graph is then encoded into the graph6 format using this
labelling. Testing whether an isomorphic copy of the new graph was already
present in the database then becomes the simple problem of checking whether the
canonically labelled graph encoded in graph6 format was already present in the
database.

It is also important that each graph is coupled with a reason for why it is considered
interesting. As already mentioned in Section~\ref{sec:comm-contr}, we therefore enforce the constraint that each newly inserted graph
is provided with a comment indicating why it was uploaded.

For each graph, \thog also provides invariant values for a selected list of invariants. 
While graphs are
submitted by users, users do not provide values for the invariants of the
uploaded graphs. Instead, these invariant values are computed by \thog on the
server (and then stored so they only have to be computed once, as many
invariants are computationally very expensive to compute). One of the reasons that the
invariant values are computed on the server rather than provided by the users, is
that it is much easier to ensure correctness and coherency of values when only
one implementation is used. As a concrete example of a coherency issue, one might
consider the diameter of a disconnected graph to be \textit{infinite} while others
might say it is \textit{undefined}. Correctness is also difficult to ensure from user-provided data. While earlier research shows that this is in principle possible in practice using
certifying algorithms~\cite{bauer2024incorporating,mcconnell2011certifying}, such practices remain,
unfortunately, uncommon at the time of writing.

Since the canonical form and the invariant values of the graphs are all computed on the server (i.e., they are not provided by the users), we only have to check if the graph was correctly formatted (no loops, undirected and within order limits). This is of course much easier than if we would have to verify the graph invariant values (or other data) inputted by the users.

\begin{figure}[h!]
  \centering
  \includegraphics[width=0.7\textwidth]{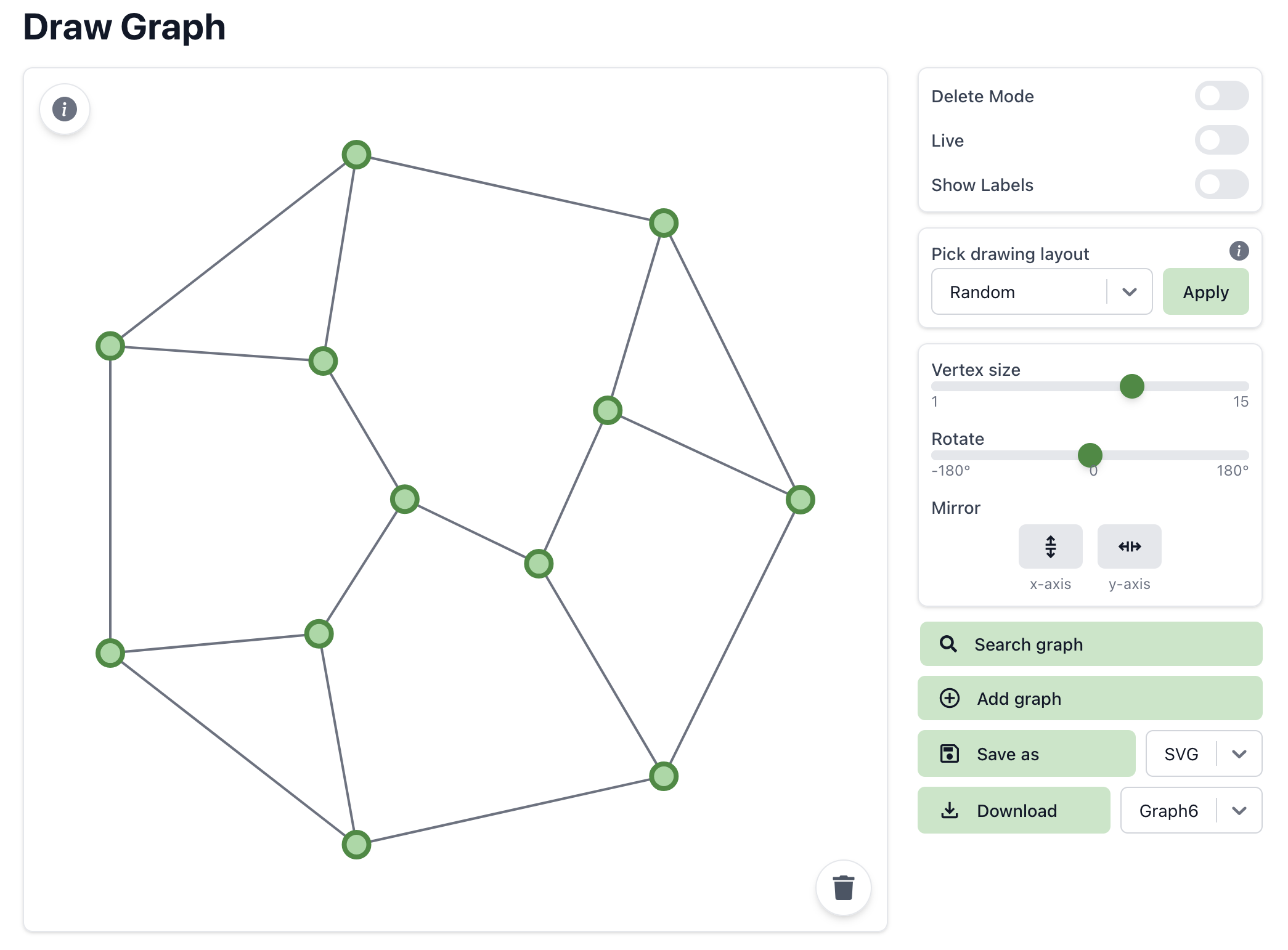}
  \caption{Overview of the graph drawing page.}
  \label{fig:draw}
\end{figure}

Users can upload a graph by either drawing it (as shown in Figure~\ref{fig:draw}) or by importing it from a file in one of the supported formats (e.g.\ graph6 or as adjacency matrix). In the latter case, a drawing of the graph is automatically generated using a spring embedding heuristic (with also the option to select other drawing heuristics as well as a heuristic specifically for planar graphs), so we can be sure that the drawing is correct. Once the graph is added, other users (or the same user) can add additional drawings for that graph by starting from the original drawing and dragging the vertices using the graph editor. However, they cannot add or remove vertices or edges, so the drawing remains correct. When the embedding is submitted, only the positions of the vertices are sent to the server, preventing an invalid drawing.

To compute the invariant values, we use well-tested software libraries -- such as \textit{nauty}~\cite{nauty} or the invariant computers from~\cite{goedgebeur2024k2} -- that are often used in research. For certain graph invariants, we have multiple independent algorithms to compute those invariants. In particular, this is the case for: chromatic index, chromatic number, circumference, domination number, edge connectivity, genus, hamiltonicity, independence number, length of longest induced cycle \& path, planarity, and vertex connectivity. In most cases these additional algorithms were added to improve the performance of the earlier invariant computers, each one smarter/faster than the previous one and sometimes also implemented in a more performant programming language (e.g.\ Java vs C++).
Ideally, we should run both algorithms and verify that they yield the same result. However, we currently do not systematically do that yet (one of the reasons being that one implementation is often a lot slower than the other and we want to avoid overloading the server, see also Section~\ref{sec:system-stability}).

Computing all invariant values on the server of course requires time limits as
many invariants are known to be NP-hard. If an invariant for a given graph
cannot be computed within the time limits (typically they are set to
approximately 5 minutes per invariant per graph), this is clearly indicated on
the page listing the graph and its variant values (cf.\
Figure~\ref{fig:petersen}) and this graph will not appear in queries concerning
one of the timed-out invariants.

 Table~\ref{tab:timeouts} shows how many graphs have a certain amount of invariants whose computation resulted in a timeout (recall that \thog currently supports 52 invariants and that the default timeout is 5 minutes).
 Table~\ref{tab:timeouts} shows that the larger numbers of timeouts are caused by a small number of graphs. Most of these graphs have large order (more than 200 vertices), but some are rather small. 
 The smallest graph having 12 timeouts has id 52206 and only 48 vertices\footnote{This graph can be accessed here: \url{https://houseofgraphs.org/graphs/52206}}. 

\begin{table}[!htbp]
\centering
\begin{tabular}{w{c}{3cm}|w{c}{3cm}}
    Number of timeouts & Number of graphs  \\ \hline
 0         & 16052   \\
 1         & 6413    \\
 2 -- 5    &  2252   \\
 6 -- 10   &  1480   \\
 11        &  302    \\
 12        & 1409    \\
 13        &  608    \\
 14        &  186    \\
 15 -- 17  &    31   \\
 $> 18$    &   0 \\
\end{tabular}
\vspace{.2cm}
\caption{Number of graphs with a given number of timed out invariant values.}
\label{tab:timeouts}
\end{table}

The timeouts are also distributed between a rather small number of invariants as seen in Table~\ref{tab:inv_timeout}, which shows the top five invariants in terms of number of timeouts. 

\begin{table}[]
    \centering
    \begin{tabular}{w{l}{5cm}|w{c}{3.5cm}}
    Invariant & Number of timeouts \\
    \hline
 Genus                           &    11513 \\ 
 Treewidth                       &     5393 \\
 Length of Longest Induced Cycle &     4715 \\
 Length of Longest Induced Path  &     4647 \\
    Hypotraceable                &     3852 \\
    \end{tabular}
    \vspace{.2cm}
    \caption{Top five invariants in terms of timeout.}
    \label{tab:inv_timeout}
\end{table}

Filling out some of this missing data could be done by selectively allocating more computational power to graphs that warrant it. However, besides our current limited budget, it is not sufficient to just increase the computational power or we would fall back into the initial problem of limited resources (albeit with a higher limit). 
Clear conditions must first be defined that determine which pairs of graphs and invariants should be assigned more computational power. Such constraints would be related to what makes the graph interesting, thus being complex to automate and calling for manual curation. 

\subsection{Relevance of the information}\label{sec:pert-inform}

Besides being correct and coherent, the information must also be relevant and of interest
for the community. As it is not feasible to store all possible
information about every graph, we restrict ourselves to information that
is useful and agreed upon by the majority of users.

The invariants stored in \thog are chosen because they are popular in the
literature or in a wide research field. 
We listen to suggestions from users and
extend the list of supported invariants if a suggestion fits this description (and if an efficient and trustworthy implementation of an algorithm to compute this invariant can be foreseen).
At the time of writing, we support 52 invariants and we believe that this includes the most commonly used invariants such as chromatic number, chromatic index, clique number, circumference, diameter, genus, treewidth, etc.

The name of graphs can be another point of friction. Indeed, graphs can bear
multiple names in different fields or subfields.
For example, the \textit{claw graph} is used to define a class of graphs that do not contain this graph as a subgraph (i.e., claw-free graphs). However, it is also the star graph on four vertices to those interested in trees. Thus, each graph can be given
a name but this name is optional (and more names can be added in the comments). Therefore, the name is not used as a unique and
deterministic identifier for graphs.  Instead, each graph is assigned a
numerical identifier generated by \thog. This identifier is unique and is used to directly link to the
graph page. This allows users to use their favourite name and still provide a
direct reference to the graph via its unique identifier without mentioning a different name.

\section{Knowledge exploration}\label{sec:exploration}

An important aspect of \thog is how to explore the knowledge. \Thog tries to
provide powerful search options while also taking limits into account to prevent abuse or
overloading the server. This is described in Section~\ref{sec:web-interface}. For the most
tech-savvy users, \thog has been integrated into some external systems via an HTTP-based
API. The specific challenges this raises are explained in Section~\ref{sec:api}.

\subsection{Web interface}\label{sec:web-interface}

One of the difficulties of such a service like House of Graphs, is to balance user freedom with safety
and stability. \Thog uses a three-component design composed of a frontend (the
web page), a backend (the business logic), and the database. The
frontend runs on the user's computer while the backend and the database execute
on \thog' server.  Part of the design of the web interface is to decide which
component should be responsible for each type of query.

\Thog stores its data in a relational database, enabling efficient SQL queries.
This already enables some simple search options for graphs, such as searching for
graphs with a given invariant value.

\begin{figure}[h!]
  \centering
  \includegraphics[width=0.8\textwidth]{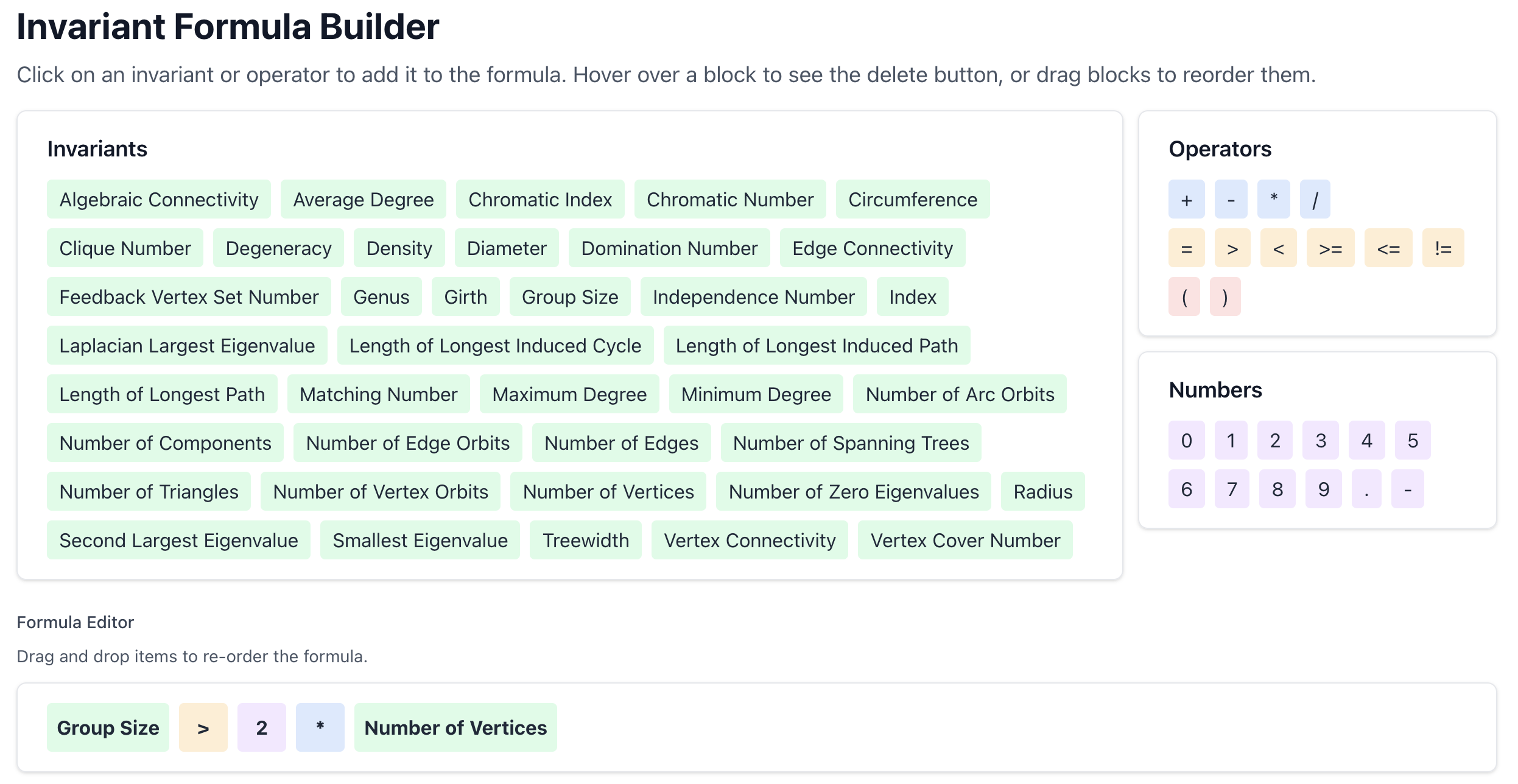}
  \caption{Overview of the search-by-formula page.}
  \label{fig:search_formula}
\end{figure}

Some queries that would be easily integrated into SQL, such as searching for a
graph whose invariants satisfy a given formula (e.g.\ to search for graphs where the order of their automorphism
group is at least as large as their number of vertices), actually require a careful User Interface
design. Indeed, one cannot expect users to know SQL syntax (nor allow users to
directly write their own SQL query for security reasons)\footnote{We realise that some users would likely be able to use SQL and that this would certainly be a powerful feature for them. However, such practice can open the door to serious security issues and is strongly discouraged by the OWASP foundation: \url{https://cheatsheetseries.owasp.org/cheatsheets/SQL_Injection_Prevention_Cheat_Sheet.html}}.
Moreover, providing a domain-specific language requires deciding on an identifier for each invariant.
Instead, \thog provides a drag-and-drop interface where users can build their
formula by piecing together blocks from different categories: invariants,
comparisons, arithmetic operations, parentheses, etc. This removes the need to
learn new syntax or naming conventions and allows users to immediately use the
tool. This is shown in Figure~\ref{fig:search_formula}.

Other search options require additional implementation in the backend. One example is
searching for a specific graph. It is possible to specify a graph in several
formats (e.g.\ in the aforementioned graph6 format or as adjacency matrix), upload it via the web interface and check if it was already present in the database. However, this
requires additional steps from the user. \Thog also provides a graph drawing
component to easily input the graph without technical knowledge or to export it in
the desired format (see Figure~\ref{fig:draw}). 
But even then, isomorphism makes it impossible to rely
solely on SQL. As explained earlier in Section~\ref{sec:corr-coher}, we use the
\emph{nauty} library~\cite{nauty} to compute a canonical labelling for every new graph. This type of query is
then ultimately transformed into an SQL query.
To search a graph from its drawing, besides the expected checks that the adjacency matrix is valid, we also compute the canonical labelling of the provided graph. This cannot be done directly by the database and thus requires an additional non-trivial computation. Once the canonical labelling is computed (via the \textit{nauty} library), the isomorphism test is reduced to a simple equality test of canonical labellings in the database using an SQL query.

The last type of search option is out of reach for a database. That is,
computations that are complex, heavy and not easily precomputed. The only such
case in \thog is \textit{subgraph search}, i.e.\ searching for graphs which contain (or do not contain) a given graph as (induced) subgraph.
This is performed using the VF2 algorithm~\cite{vf2}. This algorithm was chosen for its efficiency and the possibility to directly include it into \thog code base without additional dependencies, thus preventing an increase in maintenance work.
Still, while computing the canonical ordering of a
graph is fast enough for graphs up to 250 vertices, testing whether a graph
contains a given arbitrary subgraph, can take prohibitively long in some cases. 
This is problematic for the user who might not want to wait that long
and potentially also for other users who might face degraded performance while the
server is busy computing.

Such a situation can be solved in two ways: moving the computation to the user's
device or enforcing time limits. We opted for the latter as the first option would
require downloading the entire database on the user's device.

The time limit is applied in such a way that it still provides as much information as
possible to the user. First, it can be set by the user within a fixed range rather than
being an arbitrarily imposed fixed value, allowing some freedom to wait longer for hopefully more
results. 
Second, the search is performed incrementally from small graphs to
larger ones (for which the subgraph searches typically take longer). This allows to gather as many results as we can before the time
limit. Thus, a timeout will still return results, albeit incomplete (i.e., all graphs found up to the moment the time limit was reached). Moreover, the subgraph search can be combined with other (easier to compute) search criteria to reduce the search space, so not the entire database has to be searched.

\subsection{Programmatic exploration through the API}\label{sec:api}

In order for the frontend (the web page) to send queries to the backend (the
server), an Application Programming Interface (API) is required. \Thog uses an HTTP-based REST API. 
Such an API can also be used
programmatically without the frontend. This structure makes it possible for
users to integrate \thog into their own systems.

One such example is \textit{GraphHarverster}~\cite{deynet2024graph}. This tool is able to extract graph
drawings from a pdf file (e.g.\ a paper) and then checks via the REST API whether they were already present in \thog and then
either offers a link to the graph page or to a page where the graph can be
submitted to \thog.

Integrating \thog into code is often realised through externally developed libraries. One such example is \textit{Lean-HoG}~\cite{bauer2024incorporating}, which aims to incorporate \thog into the Lean proof assistant.

\section{Stability}\label{sec:stability}

Stability is another important requirement for an online service. This includes both system
stability and knowledge stability. The first focuses on preventing crashes and
data loss. Knowledge stability is about ensuring that the data, when
correct, is only extended and not edited. In Section~\ref{sec:system-stability}
we address the technical details to ensure system stability. Knowledge stability
is addressed in Section~\ref{sec:knowledge-stability}.

\subsection{System stability}
\label{sec:system-stability}

Preventing data loss has been extensively studied. 
To this end, we backup the database weekly via an automated script to minimise loss in case of a problem.

Some decisions are made to prevent the server from being overloaded. Graph invariants are limited to 5 minutes of
computation time and only a limited number of CPU cores is used for the invariant computations so enough cores remain available to process user queries. This is ensured by using the \textit{Slurm Workload Manager}~\cite{slurm}, i.e., a job scheduler which supports multiple queues. 
In fact, we use a \textit{multilevel feedback queue} consisting of multiple queues with different priorities and timeouts so a few long computations will not block all other (possibly much shorter) computations for a long time. See~\cite{hog2} for details.

To limit errors in invariant values, invariant computers are carefully selected and reviewed manually (see Section~\ref{sec:corr-coher} for more details). Because invariant computations are typically the most expensive computations in \thog (which luckily only have to be done once as the computed values are saved in the database), special attention is also given to the performance of the invariant computers. Indeed, while a time limit prevents the server from being overloaded, it also means that some invariant values will be missing should the limit be reached. For example, the computationally ``hardest'' invariant on \thog is the \textit{genus}, which is only known for about 60\% of the graphs in \thog as it resulted in timeouts for 40\% of the graphs. It is a special case since for the second hardest invariant, \textit{treewidth}, the timeout rate is already down to 20\% (cf.\ Table~\ref{tab:inv_timeout}).

As such, we favour algorithms implemented in performant compiled languages such as C or C++ for
hard to compute invariants, some of which go as far as to rely on bitwise
operations for maximum performance. We also prefer implementations that gained enough maturity such as
\emph{nauty}~\cite{nauty}, often implying a long history of optimisation and testing by the
community.

For long-running systems such as \thog, updates are important. Not only to add new features, but also to patch
security issues from dependencies. However, updating dependencies can quickly
result in conflicts, e.g.\ as certain methods were refactored or became deprecated. This is especially true in the context of web
applications where development is fast-paced and libraries become obsolete
quickly. While there is no clear solution, keeping the number of dependencies to
a minimum at least limits this issue.

This trade-off also impacts the choice of graph invariant computers. For
example, some state of the art invariant computers might require adding many dependencies
which, in turn, might add many hurdles when updating the system in the future.

\subsection{Knowledge stability}
\label{sec:knowledge-stability}

Even when the data is correct and coherent, data could still be modified, sometimes for subtle
reasons, and still remain correct.

A first subtle example is related to the duplicate-free storage of graphs. As
explained earlier in Section~\ref{sec:corr-coher}, a canonical ordering of the vertices of each graph is stored to allow isomorphism checking in a simple and efficient way. However, there is no standard
method to obtain this ordering. In theory, any permutation of the vertices could
be used as the canonical one. This implies that, when updating the \emph{nauty}
library~\cite{nauty}, the ordering it produces for a given graph could be different from the
one produced by the previous version. So such an update is not to be taken lightly. Therefore, before upgrading to a newer version of \textit{nauty}, we first verify that it yields the exact
same canonical ordering for every graph present in the database. This method works for
\thog because the number of interesting graphs in the database is not that large and graphs not in
the database are not impacted by the change.

Legitimate user actions can also cause data loss. Besides deleting their own
uploaded graphs or comments, users can decide to delete their account, which could cause a cascading delete. In such
situations, graphs and comments are not deleted automatically but rather their
ownership is moved to an anonymised dummy user to preserve information.

A user account can also require being disabled for reasons independent of their
actions such as security issues in cryptographic algorithms. In this case, we
prevent logging into the account but do not delete it. In order to keep using
their account, a user will be forced to change their password through an e-mail
sent to the provided address. This was indeed the case when we switched from
\textit{MD5}, which is no longer considered
collision-resistant~\cite{md5_break}, to \textit{bcrypt}~\cite{bcrypt} as
hashing scheme for the user passwords (see~\cite{hog2} for more details).

Finally, some changes can have an impact outside of \thog. Initially, the
website used the address \url{https://hog.grinvin.org}. When \thog was rebuilt
in 2021-2022, all data (including invariant values, identifiers, comments, and user accounts) was migrated and the address was changed to
\url{https://houseofgraphs.org} and also the url structure was modified. This
change may sound inconspicuous and easily solved by keeping a temporary
redirection. However, because \thog is used in a research context where papers
are published, urls in these papers can become invalid. Careful consideration is
then required to ensure backward compatibility. In particular, every previous url
must remain valid and still point or redirect to the same page. In the old House of Graphs,
accessing the Petersen graph with id 660 used the url
\url{https://hog.grinvin.org/ViewGraphInfo.action?id=660}. In the new version,
this url now is \url{https://houseofgraphs.org/graphs/660}. The new version, however, must support both schemes.

\section{Conclusion}\label{sec:conclusion}

We have presented \thog, which is a searchable database of ``interesting'' graphs. Given the amount of information stored in the database and the fact that researchers rely on the correctness and stability of the data, we discussed knowledge management in \thog and described what measures we took to make sure all data is coherent, qualitative and stable.

In the future, we hope to further extend the correctness assurances of our data by providing certificates so users can easily verify the correctness of most of the precomputed graph invariant values, along the lines of what was done in~\cite{bauer2024incorporating}. Moreover, we also aim to extend our API so other mathematical software systems such as algebra systems, proof assistants or tools like \textit{GraphHarverster}~\cite{deynet2024graph} can interact and integrate \thog into their system.

\begin{credits}
\subsubsection{\ackname} We would like to thank Katja Ber{\v{c}}i{\v{c}} and Florian Rabe for useful suggestions. Next to that, we also thank the many people who contributed to \thog by uploading new graphs, adding comments or better graph drawings or contributed in any other way. Finally, we thank the reviewers for their comments and suggestions.

Jan Goedgebeur is supported by a grant of the Research Foundation Flanders (FWO) with grant number G0AGX24N and
by Internal Funds of KU Leuven.

\subsubsection{\discintname}
The authors have no competing interests to declare that
are relevant to the content of this article.
\end{credits}

\bibliographystyle{splncs04}
\bibliography{bibliography}

\end{document}

%% file: notations.tex
\newcommand{\thog}{the House of Graphs\xspace}
\newcommand{\Thog}{The House of Graphs\xspace}
\newcommand{\gph}{GraPHedron\xspace}